\documentclass[11pt]{amsart}
\usepackage{a4wide}
\usepackage[margin=1in]{geometry}
\usepackage{amsmath,amsthm,amsfonts,amssymb,bm, mathtools, csquotes}
\usepackage{verbatim} 
\usepackage{tikz}
\usetikzlibrary{positioning,automata}
\usepackage{color}
\usepackage{graphpap}
\usepackage{epsfig}
\usepackage{psfrag}
\usepackage{graphicx}
\usepackage[all]{xy}
\usepackage{mathtools}
\usepackage{xfrac}
\usepackage{bbm}
\usepackage[T1]{fontenc}

\usepackage{hyperref}
\usepackage{mathtools}
\usepackage{microtype}
\usepackage{lmodern}
\usepackage{xcolor}
\usepackage{arydshln}
\usepackage{yfonts}
\usepackage{breqn}
\usepackage{caption}
\usepackage{geometry}

\newtheorem{teo}{Theorem}
\newtheorem{cor}{Corollary}

\newtheorem{obs}{Remark}

\def\bbP{{\mathbb P}} 
\def\B{\textbf{B}}

\def\bbN{{\mathbb N}}
\def\bbE{{\mathbb E}}

\begin{document}
	
	\title[Phase transition through embedding]{Phase transition in preferential attachment-detachment through embedding}
	
	
	\author{Michael Hinz$^1$}
	\address{$^1$ Fakult\"{a}t f\"{u}r Mathematik, Universit\"{a}t Bielefeld, Postfach 100131, 33501 Bielefeld,
		Germany}
	\email{mhinz@math.uni-bielefeld.de}
	\author{Angelica Pachon$^2$}
	\address{$^2$ Department of Mathematics, The Computational Foundry, Swansea University Bay Campus, Fabian Way,
		Swansea, SA1 8EN, United Kingdom}
	\email{a.y.pachon@swansea.ac.uk}
	
	
	\date{\today}

\begin{abstract}We study a random graph model with preferential edge attachment and detachment through the embedding into a generalized Yule model. We show that the in-degree distribution of a vertex chosen uniformly at random follows a power law in the supercritical regime but has an exponential decay in the subcritical. We provide the corresponding asymptotics. In the critical regime we observe an intermediate decay. The regimes are clearly defined in terms of parameter ranges.
\end{abstract}

\maketitle
\section{Introduction}

We study a new preferential attachment model with edge detachment, closely linked to birth-death models in evolutionary dynamics, \cite{Maruvka, Zanette}, and somewhat different from established models, \cite{DeifjenLinholm, DorogovMendes2000, Wuethal09, 2013_Vallier}. We use an embedding into a generalized Yule model, \cite{politosacerdotelanski}, to show that there are different asymptotic regimes for the tail of the in-degree distribution. Our main point is that through this embedding, the description of these different regimes is particularly neat, and our results may be seen as another illustration for the efficiency of embedding techniques, \cite{Athreya, Athreya2, Bhamidi, RSA:RSA20137}. 

Preferential attachment models are well-studied discrete-time random graph models in which new vertices are successively added and connected to an existing vertex with a probability proportional to an increasing function of its in-degree. This mechanism is well-suited to model a diverse range of phenomena, such as the concentration of resources in an market environment or the concentration of links in a communication network. The vertices represent market participants or network nodes, while the incoming edges represent the allocation of resources or directed communication channels. Such applications motivate to also include a mechanism of edge detachment: It is not atypical that the allocation of goods from one party to another stops due to changed market conditions or that communication links disappear spontaneously by a sudden loss of attention. 

Although preferential attachment models are relatively new, related models had been known long since: In \cite{Yule1925} Yule studied a continuous-time model for the creation of families and the evolution of  individuals within them. In \cite{SIMON1955} Simon studied a discrete-time model to describe the appearance of new words in a large piece of a text. For specific parameters the limit distribution of the number of individuals within a family chosen with uniform probability in the Yule model, as the number of families goes to infinity, coincides with the limit distribution of the number of occurrences of each word in the Simon model as the number of words goes to infinity. This distribution is known as the Yule-Simon distribution. In the recent article \cite{politosacerdotelanski} a generalized Yule model was proposed, where the pure birth process governing the growth of individuals within a family was replaced by a birth-death process with fertility and mortality rates $\lambda$ and $\mu$, respectively. The authors obtained the limit distribution of the size of a randomly chosen family and observed different behaviour in the supercritical regime $\lambda>\mu$, the subcritical regime $\lambda<\mu$ and the critical regime $\lambda=\mu$. Their simulations show a heavy tail power-law in the supercritical regime, a much more rapid decay in the subcritical and an intermediate decay in the critical.

The first references for preferential attachment models in the context of random graphs were \cite{342} and \cite{Barabasi1999}.
The limit degree distribution of these models as the number of vertices goes to infinity was studied rigorously in \cite{Bollobas2001}; it coincides with the Yule-Simon distribution for specific parameters.  Numerous generalizations of the model in \cite{Barabasi1999} were investigated, see for example \cite{Berger05degreedistribution,CooperFrieze, RSA:RSA20318, Deijfen09,CPC:9224998, redner,PhysRevLett.85.4629, PhysRevLett.85.4633, oliveira2005,PachonSacerdoteYang}. Links to continuous-time Markov processes were studied in \cite{Athreya, Athreya2, Bhamidi, RSA:RSA20137}. Models with edge detachment were considered in \cite{ChungLu2004, CooperFriezeVera, DeifjenLinholm,DorogovMendes2000, Wuethal09}. In \cite{ChungLu2004}
and \cite{CooperFriezeVera} it was shown that the degree distribution follows a power law whose exponent depends on the attachment and detachment probabilities. In \cite{DorogovMendes2000} it was shown that detachment changes the graph more drastically than mere attachment. In \cite{DeifjenLinholm} and \cite{Wuethal09} the authors observed a phase transition in the expected empirical degree distribution with a power law decay in the supercritical regime and an exponential decay in the subcritical. The model in \cite{Wuethal09} is very similar to special cases in \cite{ChungLu2004} and \cite{CooperFriezeVera}. Both \cite{DeifjenLinholm} and \cite{Wuethal09} gave a critical threshold, but only an incomplete description of the possible parameter regimes. A complete description was provided in \cite{2013_Vallier}, including the observation that the decay at the critical threshold is \enquote{strictly intermediate}.

Here we introduce a different model: At each time step we either add a vertex and a loop attached to it or create an edge following preferential attachment or delete an edge following preferential detachment. A precise description is given in Section \ref{TheModel}. The parameters in our model are not probabilities but rates, and the preferential attachment-detachment probabilities are obtained as the products of these rates with the percentage of oriented edges pointing to the respective vertex. By this construction the model embeds naturally into a generalized Yule model, and by this embedding it inherits the limit behaviour proved in \cite{politosacerdotelanski}. The embedding is stated in Theorem \ref{embeddingTeo} and might be of independent interest. We prove a complete description of the different asymptotic regimes for the limit in-degree distribution of a vertex chosen uniformly at random: Power law decay in the supercritical, exponential decay in the subcritical and a \enquote{strictly intermediate} decay in the critical. The parameter ranges for the different regimes admit clear interpretations in terms of evolutionary dynamics. In contrast to the aforementioned models, the exponent of the power law in our model can take any value in $(1,+\infty)$. These findings, together with results for the corresponding expected value, are stated in our main result, Theorem \ref{teo1} below.

Our model is partially motivated by two specific birth-death models, \cite{Maruvka, Zanette}. The model in \cite{Zanette} starts with one individual. At each time step a new individual is born and, moreover, an individual randomly chosen from the existing population dies with probability $\mu$. The population is divided into families, within which all individuals bear the same surname. Each newborn is assigned a new surname, so far not present in the population, with probability $\alpha$. With probability $1-\alpha$ an existing individual, chosen uniformly at random, gives their surname to the newborn. For $\mu=0$ the model reduces to Simon's model, \cite{SIMON1955}.
The model in \cite{Maruvka} describes a population of individuals (species), each of which belongs to a family (genus). It starts with a single individual.  At each time step an individual is chosen at random. This individual dies with probability $1-p$ and reproduces with probability $p$. The offspring belongs to the same family as its parent individual with probability $1-\alpha$, and \enquote{mutates} to form a new family with probability $\alpha$. 

In \cite{Deijfen10} a continuous-time network model was studied in which the vertex population evolves as a supercritical birth-death process. When a new vertex appears, an existing vertex connects to it at a rate which is a non-decreasing function of the fitness of the existing vertex. An existing vertex dies at a rate which, too, is non-decreasing in its fitness. The fitness of a vertex is its cumulative in-degree up to the present time, it takes into account current links and \enquote{ghost links} to meanwhile deceased vertices. Ghost links can model former contracts or expired communication channels. The author proved a similar phase transition between power law, exponential and intermediate regimes. For a corresponding model without ghost links they provided numerical simulations, but no rigorous calculations. Our model is based on the present in-degree, ghost links are excluded.


In Section \ref{TheModel} we provide a description of the model, in Section \ref{S:asymp} we state and discuss our main result on the limit in-degree distribution, Theorem \ref{teo1}. To set up notation and since we could not find anything comparable in the literature, we provide an inductive construction of a generalized Yule model in Section \ref{Yulemodel}. In Section \ref{S:TheEmbedding} we provide the embedding result, Theorem \ref{embeddingTeo}, and a more detailed variant of Theorem \ref{teo1}, see Theorem \ref{teo2}. To make the article largely self-contained, we recall some known facts in two Appendices.

Given functions $f,g:(0,+\infty)\to\mathbb{R}$ or $f,g:\mathbb{N}\to\mathbb{R}$, we write 
$f\sim g$ if $\lim_{t\rightarrow\infty}  \frac{f(t)}{g(t)} = 1$.
We use $\B$ to denote the Beta function,
\[\B(a,b)=\int_0^1s^{a-1}(1-s)^{b-1}\:ds,\quad a>0, b>0,\] 
and $\mathbf{U}$ to denote the confluent hypergeometric function,
\[\mathbf{U}(a,b,z)=\frac{1}{\Gamma(a)}\int_0^\infty e^{-zs}s^{a-1}(1+s)^{b-a-1}ds,\quad a>0, b>0, z>0;\]
here $\Gamma$ is the Euler gamma function. By ${}_2\mathbf{F}_1$ we denote the Gauss hypergeometric function,
\[{}_2\mathbf{F}_1(a,b,c,z)=\frac{1}{\mathbf{B}(b,c-b)}\int_0^1s^{b-1}\frac{(1-s)^{c-b-1}}{(1-sz)^a}ds,\quad a>0, b>0, c>0, 0<z<1.\]

\section{A random graph model with edge detachment}\label{TheModel}

The model we study can be described as a random graph process $(G_t^{\lambda_1,\lambda_2,\mu_2})_{t\in\mathbb{N}}$ with discrete time $t\in \mathbb{N}$ and with parameters $ \lambda_1,\lambda_2>0$ and $\mu_2\geq 0$.  At each time $t$, a realization of the process is a directed graph $G_t^{\lambda_1,\lambda_2,\mu_2}(\omega)$ with vertex set $V_t(\omega)\subset \mathbb{N}\setminus \{0\}$. The process starts at time $t=0$ as a graph $G_0^{\lambda_1,\lambda_2,\mu_2}$ consisting of the single vertex $1$ and a directed loop attached to it. The number of vertices at time $t$ is given by the random variable $|V_t|$; here $|M|$ denotes the cardinality of a finite set $M$. For any $t$ and any $i\in \mathbb{N}\setminus \{0\}$ we set $d_t^i=0$ on the event $\{|V_t|<i\}$, while on the event $\{|V_t|\geq i\}$ we define $d_t^i$ to be the in-degree of the vertex $i$ at time $t$. The transition from $G_t^{\lambda_1,\lambda_2,\mu_2}$ to $G_{t+1}^{\lambda_1,\lambda_2,\mu_2}$ is governed by the given parameters: 
The parameter $\lambda_1$ is a rate for the appearance of new vertices, $\lambda_2$ is a rate for the appearance of new edges and $\mu_2$ a rate for the removal of edges. At time $t+1$ exactly one of the following three events happens:  
\begin{enumerate} 
	\item[(1)] A new vertex and a directed loop attached to it are added to $G_t^{\lambda_1,\lambda_2,\mu_2}$. This event occurs with probability
	$|V_t| \lambda_1\big(\sum_{j=1}^{|V_t|}(\lambda_1+\lambda_2 d_t^j+\mu_2 d_t^j))^{-1}$.
	\item[(2)] A directed edge emanating from the last vertex $|V_t|$ of $G_t^{\lambda_1,\lambda_2,\mu_2}$
	and arriving at one of the vertices $1,...,|V_t|$ of $G_t^{\lambda_1,\lambda_2,\mu_2}$ is added; here vertex $1\leq i\leq |V_t|$ is chosen with  
	probability 
	$\lambda_2 d_t^i \big(\sum_{j=1}^{|V_t|}(\lambda_1+\lambda_2 d_t^j+\mu_2 d_t^j)\big)^{-1}$
	\item[(3)] A directed edge is removed from $G_t^{\lambda_1,\lambda_2,\mu_2}$; here the probability to choose an incoming edge adjacent to vertex $1\leq i\leq |V_t|$ is 
	$\mu_2 d_t^i\big(\sum_{j=1}^{|V_t|}(\lambda_1+\lambda_2 d_t^j+\mu_2 d_t^j)\big)^{-1}.$
\end{enumerate} 
%

To provide a formal description, we restrict attention to the vectors $(d_t^1,d_t^2,...d_t^{|V_t|})$ of vertex in-degrees of the random graphs $G_t^{\lambda_1,\lambda_2,\mu_2}$, $t\in\mathbb{N}$. Let 
\begin{equation}\label{E:words}
	S:=\bigcup_{n\geq 1} \mathbb{N}^n.
\end{equation}
Given $\mathbf{x}\in S$, we write 
\begin{equation}\label{E:length}
	\ell(\mathbf{x}):=n\quad \text{if}\quad \mathbf{x}\in \mathbb{N}^n,
\end{equation}
that is, if $\mathbf{x}=(x^1,x^2,...,x^n)$, and we call $\ell(\mathbf{x})$ the \emph{length} of $\mathbf{x}$. By $\delta_{(1)}$ we denote the point mass probability measure on $S$ at the sequence $(1)$ of length one and with single element $1$. We use the convenient shortcut notation 
\begin{equation}\label{E:norm}
	\|\mathbf{x}\|:=\sum_{j=1}^{\ell(\mathbf{x})} x^j,\quad \mathbf{x}\in S.
\end{equation}
Given $\mathbf{x},\mathbf{y}\in S$ let 
\[p_{\mathbf{x},\mathbf{y}}:=\begin{cases} 
	\frac{x^i\lambda_2}{(\lambda_2+\mu_2)\|\mathbf{x}\|+\ell(\mathbf{x})\lambda_1} &\text{if $\mathbf{y}=(x^1,x^2,...,x^{i-1},x^i+1,x^{i+1},...,x^{\ell(\mathbf{x})}),\quad 1\leq i\leq \ell(\mathbf{x})$,}\\
	\frac{x^i\mu_2}{(\lambda_2+\mu_2)\|\mathbf{x}\|+\ell(\mathbf{x})\lambda_1} &\text{if $\mathbf{y}=(x^1,x^2,...,x^{i-1},x^i-1,x^{i+1},...,x^{\ell(\mathbf{x})}),\quad 1\leq i\leq \ell(\mathbf{x})$,}\\
	\frac{\ell(\mathbf{x})\lambda_1}{(\lambda_2+\mu_2)\|\mathbf{x}\|+\ell(\mathbf{x})\lambda_1} &\text{if $\mathbf{y}=(x^1,x^2,...,x^{\ell(\mathbf{x})},1),\quad 1\leq i\leq \ell(\mathbf{x})$,}
\end{cases}\]

Now let 
\begin{equation}\label{MarkovChain}
	\mathbf{d}=(\mathbf{d}_t)_{t\in\mathbb{N}}\quad \text{with}\quad \mathbf{d}_t=(d_t^1,d_t^2,...,d_t^{\ell(\mathbf{d}_t)}), 
	\ t\in\mathbb{N},
\end{equation}
be the uniquely determined discrete-time Markov chain with state space $S$, transition probabilities $(p_{\mathbf{x},\mathbf{y}})_{\mathbf{x},\mathbf{y}\in S}$  and initial distribution $\delta_{(1)}$. Ionescu-Tulcea's theorem gua\-rantees the existence of a probability space $(\Omega,\mathcal{F},\mathbb{P})$ such that $\mathbb{P}(\mathbf{d}_0=(1))=1$
and 
$\mathbb{P}(\mathbf{d}_{t+1}=\mathbf{y}|\mathbf{d}_t=\mathbf{x})=p_{\mathbf{x},\mathbf{y}}$
for all $\mathbf{x},\mathbf{y}\in S$ and $t\in\mathbb{N}$. We write $\bbE$ and $\mathbb{V}$ for the expectation and variance associated with $\bbP$, respectively.

\section{Limit in-degree distributions and tail behaviour}\label{S:asymp}

For any $t\in \mathbb{N}$ let
\begin{equation}\label{EmpiricalDistj}
	\pi_t^j=\frac{\sum_{i=1}^{\ell(\mathbf{d}_t)} \mathbbm{1}_{\{d_t^i=j\}}}{\ell(\mathbf{d}_t)},\quad j\in\mathbb{N}.
\end{equation}
Since $\ell(\mathbf{d}_t)$ equals the number $|V_t|$ of vertices in $G_t^{\lambda_1,\lambda_2,\mu_2}$, the random variable $\pi_t^j$ may be viewed as the number of vertices in $G_t^{\lambda_1,\lambda_2,\mu_2}$ having in-degree $j$. Let $d(V_t)$ be the random variable on $(\Omega,\mathcal{F},\mathbb{P})$ fully determined by the set of conditional distributions 
\begin{equation}\label{E:conditional}
	\bbP(d(V_t)=j\mid \ell(\mathbf{d}_t)=n)=\frac{1}{n}\sum_{i=1}^n\bbP(d_t^i=j \mid \ell(\mathbf{d}_t)=n),\quad n\in\mathbb{N}\setminus \{0\}.
\end{equation}
It may be viewed as the in-degree of a vertex chosen uniformly at random in $G_t^{\lambda_1,\lambda_2,\mu_2}$. Conditioning on 
$\ell(\mathbf{d}_t)$ and using (\ref{E:conditional}), we obtain  
\begin{multline}
	\bbE[\pi_t^j]
	=\sum_{n=1}\bbE\left[\left(\frac{\sum_{i=1}^{\ell(\mathbf{d}_t)}\mathbbm{1}_{\{d_t^i=j\}}}{\ell(\mathbf{d}_t)}\right)\Big| \ell(\mathbf{d}_t)=n\right]\bbP(\ell(\mathbf{d}_t)=n)\notag\\
	=\sum_{n=1}\frac{1}{n}\Big(\sum_{i=1}^n\bbP(d_t^i=j\mid \ell(\mathbf{d}_t)=n)\Big)\bbP(\ell(\mathbf{d}_t)=n)
	=\bbP(d(V_t)=j).
\end{multline}		
That is, at any fixed time $t\in\mathbb{N}$, the expectation under $\mathbb{P}$ of the empirical in-degree distribution $(\pi_t^j)_{j\geq 0}$ and the in-degree distribution $\bbP(d(V_t)=\cdot)$ under $\mathbb{P}$ of a vertex chosen uniformly at random coincide. 

Our main statement is about the limit of these in-degree distributions as $t\to\infty$. It says that a limit distribution exists and follows different asymptotic regimes as $j\to\infty$.

\begin{teo}\label{teo1}
	Let $\lambda_1,\lambda_2>0$ and $\mu_2\geq 0$, let $\mathbf{d}$ be as in (\ref{MarkovChain}) and let $d(V_t)$ be the random variable determined by (\ref{E:conditional}). Then the limits
	\begin{align}\label{ExpEmpiricalDist}
		p_j:=\lim_{t\rightarrow\infty}\bbP(d(V_t)=j),\quad j\in\mathbb{N},
	\end{align}
	exist. As $j\to\infty$, the limit distribution $(p_j)_{j\geq 0}$ shows the following asymptotic behaviour:
	\begin{enumerate}
		\item[(i)] If $\mu_2=0$, then 
		\[p_j\sim c_1\:j^{-\Big(1+\frac{\lambda_1}{\lambda_2}\Big)}\quad \text{with}\quad  c_1=\frac{\lambda_1}{\lambda_2}\Gamma\Big(1+\frac{\lambda_1}{\lambda_2}\Big).\]
		\item[(ii)] If $\mu_2>0$ and $\lambda_2>\mu_2$, then \[p_j \sim c_2 j^{-\left(1+\frac{\lambda_1}{\lambda_2-\mu_2}\right)}\quad \text{with}\quad c_2=\frac{\lambda_1}{\lambda_2}
		\Big(\frac{\lambda_2}{\lambda_2-\mu_2}\Big)^{\frac{\lambda_1}{\lambda_2-\mu_2}}\Gamma\Big(1+\frac{\lambda_1}{\lambda_2-\mu_2}\Big).\]
		\item[(iii)] If $\mu_2>0$ and $\lambda_2<\mu_2$, then 
		\[p_j\sim  c_3\:\left(\frac{\lambda_2}{\mu_2}\right)^{j+1} j^{-\left(1+\frac{\lambda_1}{\mu_2-\lambda_2}\right)}\quad \text{with}\quad c_3=\frac{\lambda_1}{\mu_2}\Big(\frac{\mu_2}{\mu_2-\lambda_2}\Big)^{-\frac{\lambda_1}{\mu_2-\lambda_2}}\Gamma\Big(1+\frac{\lambda_1}{\mu_2-\lambda_2}\Big).\]
		\item[(iv)] If $\mu_2>0$ and $\lambda_2=\mu_2$, then $j^mp(j)\to 0$ for any $m\in \mathbb{N}$ but $e^{\varepsilon j} p(j)\to +\infty$ for any $\varepsilon>0$.
	\end{enumerate}
	For the expectation of a random variable $d(V_\infty)$ with distribution $(p_j)_{j\geq 0}$ we have
	\begin{equation}\label{Expectation}
		\mathbb{E}[d(V_\infty)]=
		\begin{cases}
			\frac{\lambda_1}{\lambda_1-(\lambda_2-\mu_2)}, & \text{ if }  \lambda_2<\mu_2    \text{ or }       0<\lambda_2-\mu_2<\lambda_1  \\
			1, & \text{ if } \lambda_2=\mu_2\\
			\infty, &\text{ otherwise}.
		\end{cases}
	\end{equation}
\end{teo}


Theorem \ref{teo1} follows from Theorem \ref{teo2} stated in Section \ref{S:TheEmbedding} and known asymptotics for $\mathbf{B}$, ${}_2\mathbf{F}_1$ and $\mathbf{U}$. The proof of Theorem \ref{teo2} is based on an embedding of the Markov chain $\mathbf{d}$ into a generalized Yule model and on the asymptotics forthe latter proved in \cite{politosacerdotelanski}.

\begin{obs}
	In regime (i) \emph{no edge detachment} happens and $(p_j)_{j\geq 0}$ is a Yule-Simon distribution. In the \emph{supercritical regime} (ii) edge detachment happens, but attachment outweighs it, $\lambda_2>\mu_2$. In this case we observe a power law behaviour of $(p_j)_{j\geq 0}$. In the \emph{subcritical regime} (iii) detachment outweighs attachment and the tail of $(p_j)_{j\geq 0}$ decays exponentially. In the \emph{critical regime} (iv) we observe an intermediate decay, strictly faster than polynomial and strictly slower than exponential.
	
	If the detachment rate equals the positive attachment rate, $\mu_2=\lambda_2$, then the \enquote{expected degree at infinity} $\mathbb{E}[d(V_\infty)]$ is one. If detachment dominates attachment, $\lambda_2<\mu_2$, or the rate at which new vertices appear dominates a positive bias towards attachment, $0<\lambda_2 -\mu_2<\lambda_1$, then $\mathbb{E}[d(V_\infty)]$ is smaller than one. In the first case this is linked to the exponential decay in the subcritical regime (iii), while the second case concerns the supercritical regime, but $jp_j\sim j^{-\lambda_1/(\lambda_2-\mu_2)}$ is summable. The transition between these cases is continuous. If no detachment happens, $\mu_2=0$, and $\lambda_2<\lambda_1$, then again $\mathbb{E}[d(V_\infty)]$ is less than one, and in particular, $jp_j\sim j^{-\lambda_1/\lambda_2}$ in (i) is summable. If $\mu_2=0$ and $\lambda_2\geq \lambda_1$, or if attachment dominates a positive detachment, $\lambda_2>\mu_2>0$, and the rate at which new vertices appear is rather small, $\lambda_1\leq \lambda_2-\mu_2$, then too many edges are produced to have a finite expectation. In the first case the summability of $jp_j$ is lost in (i), in the second case it is lost in (ii).
\end{obs}

\section{A generalized Yule model}\label{Yulemodel}

Recall that a \emph{birth-death process with birth rate $\lambda>0$ and death rate $\mu\geq 0$}, is a continuous-time Markov chain $(Z_u)_{u\geq0}$  with state space $\bbN$ and infinitesimal generator $A:=\left((a_{ij})\right)$, where $a_{ii}=-i(\lambda+\mu)$, $a_{ij}= i\lambda$ if $j=i+1$, $a_{ij}= i\mu$ if $j=i-1$, and $a_{ij}=0$ if $j\neq i,i+1,i-1$. 
If $\mu=0$, then it is also called a \emph{pure birth process with birth rate $\lambda>0$}. We say that a birth-death process $(Z_u)_{u\geq0}$ \emph{starts at one} if $Z_0=1$.


Recall that for a birth-death process $(Z_u)_{u\geq0}$ the transition from a state $k>0$ can go either to state $k+1$ or to state $k - 1$, at exponential rates $k\lambda$ or $k\mu$, respectively. The waiting time between the last transition to $k>0$ and the transition from $k$ to $k+1$ or $k-1$ is independent of the history up to this last transition and exponentially distributed with parameter $k(\lambda+\mu)$.

A \emph{generalized Yule model with detachment}, \cite{politosacerdotelanski}, and \emph{with parameters $\lambda_1,\lambda_2>0$ and $\mu_2\geq 0$} is a combination of a pure birth process and a sequence of birth-death processes. It models the evolution of a population divided into households (or families). The pure-birth process describes the growth of the number of households, while the birth-death processes describes the evolution of the numbers of individuals within the particular households. Informally speaking, the model works as follows:
\begin{enumerate}
	\item Households appear according to a pure birth process with parameter  $\lambda_1$  and starting at one.
	\item Each time a new household appears, it starts a copy of  a birth-death process  with parameters $ \lambda_2$ and $\mu_2$, and starting at one.
\end{enumerate}  
This gives rise to a sequence of birth-death processes, which grow independently from each other and from the birth process modeling the appearance of households. The case $\mu_2=0$ reproduces the classical Yule model, \cite{Yule1925}.

Let $S$ be as in (\ref{E:words}). Below we construct a generalized Yule model with detachment as  a continuous-time $S$-valued process $\hat Z=(\hat{Z}_u)_{u\geq 0}$,
\begin{equation}\label{E:Z}
	\hat{Z}_u=(\hat{Z}^1_u,\hat{Z}_u^2,...,\hat{Z}_u^{\ell(\hat{Z}_u)}),\quad u\geq 0;
\end{equation}
here we use notation (\ref{E:length}). The random vector $\hat{Z}_u$ describes a population at time $u$, which at that time
consists of $\ell(\hat{Z}_u)$ separate households of sizes $\hat{Z}_u^j$, $j\geq 1$. As the process $\hat Z$ evolves, new households appear and, independently of the formation of these new households, individuals within each household are born or die.  Let $\hat{Z}^j=(\hat{Z}^j_u)_{u\geq 0}$, $j\geq 1$, be a sequence of independent copies of a birth-death process $Z=(Z_u)_{u\geq 0}$ with parameters $\lambda_2>0$ and $\mu_2\geq 0$.  At time $u=0$  the process $\hat Z$ starts with a household formed by a single individual, that is, $\hat{Z}_0=(\hat{Z}_0^1)=(1)$, and we have $\ell(\hat Z_0)=1$.  As time passes, three different events may occur: Firstly, a new individual may be born into the household counted by $\hat{Z}^1$. Secondly, an individual in this household may pass away.
The third possibility is that a new household, then counted by  $\hat{Z}^2$, is formed. This happens according to a pure birth process $\ell(\hat Z)$ with rate $\lambda_1>0$ and independent of $\hat{Z}^1$. Whatever event happens first, we label the random time it happens $\hat{\sigma}_1$. Now further births and deaths may happen in $\hat{Z}^1$, and if it was already started, also in $\hat{Z}^2$.
In addition, a new household with a single individual may be formed,  counted by another independent copy of $Z$. 
We write $\hat{\sigma}_2$ for the time that passes after $\hat{\tau}_1:=\hat{\sigma}_1$ until another birth, death or independent formation of a new household occurs, and we set $\hat{\tau}_2:=\hat{\tau}_1+\hat{\sigma}_2$. We now continue in a similar manner. This creates a sequence $(\hat{\sigma}_k)_{k\in \mathbb{N}}$ of random times with partial sums $\hat{\tau}_t=\sum_{k=1}^t \hat{\sigma}_k$. To the random times $\hat{\tau}_t$ we refer as \emph{census times}, they form an increasing sequence $(\hat{\tau}_t)_{t\in \mathbb{N}}$. We write $(\hat{\tau}^\ast_t)_{t\geq 1}$ for the subsequence of those census times at which a new houshold is formed, and we call them \emph{formation times}. The households appear according to the pure birth process $\ell(\hat Z)$ and the formation times are exactly the successive jump times of this process. Consequently $\hat{\tau}^\ast_1,\hat{\tau}^\ast_2-\hat{\tau}^\ast_1,...,\hat{\tau}^\ast_{t}-\hat{\tau}^\ast_{t-1},... $ are independent exponentially distributed random variables with parameters $\lambda_1,2\lambda_1,...,t\lambda_1,...$, respectively. 

To put this construction of $\hat{Z}$ into formulas, let $(\Omega',\mathcal{F}',\mathbb{P}')$ be the probability space over which the birth-death process $Z=(Z_u)_{u\geq 0}$ is defined and such that $\mathbb{P}'(Z_0=1)=1$. We consider the space $(\hat{\Omega},\hat{\mathcal{F}},\hat{\mathbb{P}})$, where
\[\hat{\Omega}:=\bigtimes_{j=1}^\infty \Omega'\times\bigtimes_{k=1}^\infty [0,+\infty),\quad \hat{\mathcal{F}}:=\bigotimes_{j=1}^\infty\mathcal{F}'\otimes\bigotimes_{k=1}^\infty \mathcal{B}([0,+\infty))\quad \text{and}\quad \hat{\mathbb{P}}:=\bigotimes_{j=1}^\infty\mathbb{P}' \otimes \bigotimes_{k=1}^\infty \mathcal{E}_{k\lambda_1};\]
here  $\mathcal{E}_{k\lambda_1}$ denotes the exponential distribution with parameter $k\lambda_1>0$. We use the notation $(\omega,\pi)=((\omega_1,\omega_2,...),(\pi_1,\pi_2,...))$ 
with $\omega_j\in \Omega'$ and $\pi_k\in [0,+\infty)$ for the elements of $\hat{\Omega}$. We now define $\hat{Z}$ inductively. 

Let $\hat{Z}^1=(\hat{Z}^1_u)_{u\geq 0}$ be the process over $(\hat{\Omega},\hat{\mathcal{F}},\hat{\mathbb{P}})$, defined by 
$\hat{Z}_u^1((\omega,\pi)):=Z_u(\omega_1), u\geq 0.$
It is convenient to set $\hat{\sigma}_0:=0$, $\hat{\tau}_0:=0$ and $\hat{\tau}^\ast_0:=0$. 

We define the random times
\[\sigma_1^+(\hat{Z}^1):=\inf\{u>0:\ \hat{Z}_u^1=\hat{Z}_0^1+1\},\quad  \sigma_1^-(\hat{Z}^1):=\inf\{u>0:\ \hat{Z}_u^1=\hat{Z}_0^1-1\}\] 
and 
$\sigma_1^\ast(\hat{Z})((\omega,\pi)):=\pi_{\ell(\hat{Z}_0)}=\pi_1,$
note that under $\hat{\mathbb{P}}$ the time $\sigma_1^\ast(\hat{Z})$ is independent of $\hat{Z}^1$ and exponentially distributed with parameter $\lambda_1$. 
We then set 
\[\hat{\sigma}_1:=\min\{\sigma_1^+(\hat{Z}^1),\sigma_1^-(\hat{Z}^1),\sigma_1^\ast(\hat{Z})\},\]
$\hat{\tau}_1:=\hat{\sigma}_1$ and define $\hat{Z}$ up to time $\hat{\tau}_1$ by
\[\hat{Z}_u:=(\hat{Z}^1_u),\quad \hat{\tau}_0\leq u<\hat{\tau}_1.\]

If $\hat{\sigma}_1=\sigma_1^+(\hat{Z}^1)$ or $\hat{\sigma}_1=\sigma_1^-(\hat{Z}^1)$, then we set 
$\check{Z}_u^{(1)}:=(\hat{Z}_u^1), u\geq \hat{\tau}_1.$
If  $\hat{\sigma}_1=\sigma_1^\ast(\hat{Z})$, then we set $\hat{\tau}_1^\ast:=\hat{\tau}_1$ and 
$\check{Z}_u^{(1)}:=(\hat{Z}_u^1,\hat{Z}^2_u), u\geq \hat{\tau}_1,$
where $\hat{Z}^2_u((\omega,\pi)):=Z_{u-\hat{\tau}_1^\ast}(\omega_2)$, $u\geq \hat{\tau}_1^\ast$. 
The process $(\check{Z}_u^{(1)})_{u\geq \hat{\tau}_1}$ is a \enquote{preliminary future version} of $\hat{Z}$ starting at time $\hat{\tau}_1$.

Now suppose that $t\in\mathbb{N}\setminus\{0\}$, the random time $\hat{\tau}_t$ has been determined and $(\check{Z}_u^{(t)})_{u\geq \hat{\tau}_t}$ has been defined. The vector $\check{Z}_{\hat{\tau}_t}^{(t)}$ has $\ell(\check{Z}_{\hat{\tau}_t}^{(t)})$ components,
$\check{Z}_{\hat{\tau}_t}^{(t)}=(\hat{Z}^1_{\hat{\tau}_t}, ..., \hat{Z}^{\ell(\check{Z}_{\hat{\tau}_t}^{(t)})}_{\hat{\tau}_t}).$

For any $j=1,...,\ell(\check{Z}_{\hat{\tau}_t}^{(t)})$, define the random times
\[\hat{\sigma}_{t+1}^+(\hat{Z}^j):=\inf\{u>0:\ \hat{Z}_{\hat{\tau}_t+u}^j=\hat{Z}_{\hat{\tau}_t}^j+1\},\quad  \hat{\sigma}_{t+1}^-(\hat{Z}^j):=\inf\{u>0:\ \hat{Z}_{\hat{\tau}_t+u}^j=\hat{Z}_{\hat{\tau}_t}^j-1\}\]
and 
$\sigma^\ast_{t+1}(\hat{Z})((\omega,\pi)):=\pi_{\ell(\check{Z}_{\check{\tau}_t}^{(t)})}.$
Under $\hat{\mathbb{P}}$ the time $\sigma^\ast_{t+1}(\hat{Z})$ is independent of $(\check{Z}_u^{(t)})_{u\geq \hat{\tau}_t}$ and exponentially distributed with parameter $\ell(\check{Z}_{\hat{\tau}_t}^{(t)})\lambda_1$. 
We set
\[\hat{\sigma}_{t+1}:=\min\{\sigma_{t+1}^+(\hat{Z}^1), ..., \sigma_{t+1}^+(\hat{Z}^{\ell(\check{Z}_{\hat{\tau}_t}^{(t)})}),\sigma_{t+1}^-(\hat{Z}^1), ..., \sigma_{t+1}^-(\hat{Z}^{\ell(\check{Z}_{\hat{\tau}_t}^{(t)})}), \sigma_{t+1}^\ast(\hat{Z})\},\] 
$\hat{\tau}_{t+1}:=\sum_{k=1}^{t+1}\hat{\sigma}_k$ and define $\hat{Z}$ between time $\hat{\tau}_t$ and $\hat{\tau}_{t+1}$ by
\[\hat{Z}_u:=(\hat{Z}^1_u,...,\hat{Z}^{\ell(\check{Z}_{\hat{\tau}_t}^{(t)})}_u),\quad \hat{\tau}_t\leq u<\hat{\tau}_{t+1}.\]

If for some $j=1,...,\ell(\check{Z}_{\hat{\tau}_t}^{(t)})$ we have 
$\hat{\sigma}_{t+1}=\sigma_{t+1}^+(\hat{Z}^j)$ or $\hat{\sigma}_{t+1}=\sigma_{t+1}^-(\hat{Z}^j)$, then we set\\ 
$\check{Z}_u:=(\hat{Z}^1_u,...,\hat{Z}^{\ell(\check{Z}_{\hat{\tau}_t}^{(t)})}_u), u\geq \hat{\tau}_{t+1}.$
If instead $\hat{\sigma}_{t+1}=\sigma_{t+1}^\ast(\hat{Z})$, then we set $\hat{\tau}_{t+1}^\ast:=\hat{\tau}_{t+1}$ and\\
$\check{Z}_u:=(\hat{Z}^1_u,...,\hat{Z}^{\ell(\check{Z}_{\hat{\tau}_t}^{(t)})}_u, \hat{Z}^{\ell(\check{Z}_{\hat{\tau}_t}^{(t)})+1}_u), u\geq \hat{\tau}_{t+1},$
where $\hat{Z}^{\ell(\check{Z}_{\hat{\tau}_t}^{(t)})+1}_u((\omega,\pi)):=Z_{u-\hat{\tau}_{t+1}^\ast}(\omega_{t+2})$, $u\geq \hat{\tau}_{t+1}^\ast$.

Proceeding inductively we now obtain the desired process $\hat Z=(\hat{Z}_u)_{u\geq 0}$.

\begin{obs}\label{R:strictmon}
	It is immediate that $\hat{\mathbb{P}}(\hat{\tau}_t\leq \hat{\tau}_{t+1}:\ t\in \mathbb{N})=1$. Since the sequence of jump times $(\hat{\tau}^\ast_t)_{t\in \mathbb{N}}$ of the pure birth process $\ell(\hat{Z})$ is a subsequence of $(\hat{\tau}_t)_{t\in\mathbb{N}}$, we also have $\hat{\mathbb{P}}(\lim_{t\to\infty} \hat{\tau}_t=+\infty)=\hat{\mathbb{P}}(\lim_{t\to\infty} \hat{\tau}_t^\ast=+\infty)=1$.
\end{obs}

\section{An embedding and its consequences}\label{S:TheEmbedding}

We observe the following embedding of $\mathbf{d}$ into $\hat{Z}$.

\begin{teo}\label{embeddingTeo}
	Let $\lambda_1,\lambda_2>0$ and $\mu_2\geq 0$. Let $\mathbf{d}=(\mathbf{d}_t)_{t\in\mathbb{N}}$ be the discrete-time Markov chain in (\ref{MarkovChain}), let $\hat Z=(\hat{Z}_u)_{u\geq 0}$ the continuous-time $S$-valued process
	in (\ref{E:Z}) and let $(\hat{\tau}_t)_{t\in\mathbb{N}}$ be the associated sequence of census times. Then the processes $(\hat{Z}_{\hat{\tau}_t})_{t\in\mathbb{N}}$ and $(\mathbf{d}_t)_{t\in\mathbb{N}}$ have the same distribution. 
\end{teo}
The arguments for the proof are standard, see for instance \cite[Theorem 3.3]{Athreya}. We provide some comments for the convenience of the reader. 

Recall that if $Y_1,..., Y_k$ are independent exponential random
variables with parameters  $\delta_1,..., \delta_k$, then $\min(Y_1, Y_2,... , Y_k)$ is exponentially
distributed with parameter $\delta_1+...+\delta_k$ and the probability that $\min(Y_1, Y_2,... , Y_k)=Y_i$ is $\delta_i(\sum_{j=1}^k \delta_j)^{-1}$.

\begin{proof}
	Suppose that $n\in \mathbb{N}\setminus\{0\}$, ${\bf{x}}=(x_1,x_2,...,x_n)\in S$ and $ \hat{Z}_{\hat{\tau}_{t}}=\mathbf{x}$.
	Stepping from $\hat{\tau}_{t}$ to $\hat{\tau}_{t+1}$, a birth of an individual in household $\hat{Z}^i$ occurs with probability 
	\[\hat{\mathbb{P}}\big(\hat{\sigma}_{t+1}=\sigma_{t+1}^+(\hat{Z}^i)|\ \hat{Z}_{\hat{\tau}_{t}}=\mathbf{x}\big)=\frac{x_i\lambda_2}{\sum_{j=1}^n(\lambda_2+\mu_2)x_j + n\lambda_1},\]
	a death of an individual in household $\hat{Z}^i$ with probability
	\[\hat{\mathbb{P}}\big(\hat{\sigma}_{t+1}=\sigma_{t+1}^-(\hat{Z}^i)|\ \hat{Z}_{\hat{\tau}_{t}}=\mathbf{x}\big)=\frac{x_i\mu_2 }{\sum_{j=1}^n(\lambda_2+\mu_2)x_j + n\lambda_1 }\]
	and the formation of a new household with probability
	\[\hat{\mathbb{P}}\big(\hat{\sigma}_{t+1}=\sigma_{t+1}^\ast(\hat{Z})|\ \hat{Z}_{\hat{\tau}_{t}}=\mathbf{x}\big)=\frac{n\lambda_1 }{\sum_{j=1}^n(\lambda_2+\mu_2)x_j + n\lambda_1 }.\]
	A comparison shows that these are exactly the transition probabilities of $\mathbf{d}$ in (\ref{MarkovChain}). In particular, 
	the process $(\hat{Z}_{\hat{\tau}_t})_{t\in\mathbb{N}}$ is a discrete-time Markov chain with state space $S$. Since also $\mathbf{d}_0=(1)=\hat{Z}_{\hat{\tau}_0}$ by construction, the Markov chains have the same distribution.
\end{proof}

Let $\hat Z=(\hat{Z}_u)_{u\geq 0}$ the continuous-time $S$-valued process $\hat Z=(\hat{Z}_u)_{u\geq 0}$
in (\ref{E:Z}) with parameters $\lambda_1$, $\lambda_2$ and $\mu_2$. For any $u\geq 0$, let $s(\hat{Z}_u)$ be the random variable on $(\hat{\Omega},\hat{\mathbb{P}},\hat{\mathcal{F}})$ fully determined by the set of conditional distributions 
\begin{equation}\label{E:size}
	\hat{\bbP}(s(\hat{Z}_u)=j\mid \ell(\hat{Z}_u)=n)=\frac{1}{n}\sum_{i=1}^n\hat{\bbP}(\hat{Z}_u^i=j \mid \ell(\hat{Z}_u)=n),\quad n\in\mathbb{N}\setminus\{0\}.
\end{equation}
It may be interpreted as the size  of a household  chosen uniformly at random at time $u$. The following is immediate from Theorem \ref{embeddingTeo} and Remark \ref{R:strictmon}.

\begin{cor}\label{C:asymp}
	Let $\lambda_1,\lambda_2>0$ and $\mu_2\geq 0$. Given $t\in \mathbb{N}$ and $u\geq 0$, let $d(V_t)$ and $s(\hat{Z}_u)$ be the random variables determined by (\ref{E:conditional}) and (\ref{E:size}), respectively. Then we have $\mathbb{P}(d(V_t)= j)=\hat{\mathbb{P}}(s(\hat{Z}_{\hat\tau_t})=j)$
	for any $t\in\mathbb{N}$ and $j\in\mathbb{N}$, where $(\hat{\tau}_t)_{t\in\mathbb{N}}$ is the sequence of census times associated with $\hat{Z}$. In particular, 
	\begin{equation}\label{uniformchoice}
		\lim_{t\rightarrow\infty}\bbP(d(V_t)= j)=\lim_{u\rightarrow\infty}\hat{\bbP}(s(\hat{Z}_u)= j),\quad j\in\mathbb{N}.
	\end{equation}
\end{cor}

The following Theorem \ref{teo2} is a straightforward consequence of the embedding in Theorem \ref{embeddingTeo} and Corollary \ref{C:asymp} together with 
\cite{Yule1925} and the results in [\cite{politosacerdotelanski} Sections 2 and 3]. For the convenience of the reader we recall known key facts from \cite{Yule1925} and \cite{politosacerdotelanski}
in Appendix \ref{S:key}, together with brief comments on (\ref{Expectation}).

\begin{teo}\label{teo2}
	Let $\lambda_1,\lambda_2>0$ and $\mu_2\geq 0$, let $\mathbf{d}$ be as in (\ref{MarkovChain}) and let $d(V_t)$ be the random variable determined by (\ref{E:conditional}). Then the limits $p_j$ in (\ref{ExpEmpiricalDist}) exist and are as follows:
	\begin{enumerate}
		\item[(i)] If $\mu_2=0$, then
		\begin{equation}   \label{nuZero}
			p_0=0\quad\text{and}\quad 
			p_j=\frac{\lambda_1}{\lambda_2} \B\Big(j,1+\frac{\lambda_1}{\lambda_2}\Big) 
			,\quad j\in \mathbb{N}\setminus \{0\}. 
		\end{equation} 
		\item[(ii)] If $\mu_2>0$ and $\lambda_2>\mu_2$, then 
		\begin{equation}\label{supercriticalp0}
			p_0=\frac{\mu_2\lambda_1}{\lambda_2(\lambda_2-\mu_2)}\mathbf{B}\Big(2,\frac{\lambda_1}{\lambda_2-\mu_2}\Big){}_2\mathbf{F}_1\Big(1,\frac{\lambda_1}{\lambda_2-\mu_2},2+\frac{\lambda_1}{\lambda_2-\mu_2},\frac{\mu_2}{\lambda_2}\Big), \text{ and }  
		\end{equation} 
		\begin{align}\label{supercriticalpj}
			p_j=\frac{\lambda_1(\lambda_2-\mu_2)}{\lambda_2^2}\mathbf{B}\Big(j,1+\frac{\lambda_1}{\lambda_2-\mu_2}\Big){}_2\mathbf{F}_1\Big(j+1,1+\frac{\lambda_1}{\lambda_2-\mu_2},j+1+\frac{\lambda_1}{\lambda_2-\mu_2},\frac{\mu_2}{\lambda_2}\Big),\nonumber\\j\in \mathbb{N}\setminus\{0\}  
		\end{align}
		\item[(iii)] If $\mu_2>0$ and $\lambda_2<\mu_2$, then
		\begin{equation}\label{subcriticalp0}
			p_0=\frac{\lambda_1}{\mu_2-\lambda_2}\mathbf{B}\Big(2,\frac{\lambda_1}{\mu_2-\lambda_2}\Big){}_2\mathbf{F}_1\Big(1,\frac{\lambda_1}{\mu_2-\lambda_2},2+\frac{\lambda_1}{\mu_2-\lambda_2},\frac{\lambda_2}{\mu_2}\Big), \text{ and }
		\end{equation}
		\begin{multline}\label{subcriticalpj}
			p_j=\left(\frac{\lambda_2}{\mu_2}\right)^{j+1}\frac{\lambda_1(\mu_2-\lambda_2)}{\mu_2^2}\mathbf{B}\Big(j,1+\frac{\lambda_1}{\mu_2-\lambda_2}\Big)\\{}_2\mathbf{F}_1\Big(j+1,1+\frac{\lambda_1}{\mu_2-\lambda_2},j+1+\frac{\lambda_1}{\mu_2-\lambda_2},\frac{\lambda_2}{\mu_2}\Big),  j\in \mathbb{N}\setminus \{0\}
		\end{multline} 
		\item[(iv)] If $\mu_2>0$ and $\lambda_2=\mu_2$, then 
		\begin{equation}\label{criticalp0}   
			p_0=\mathbf{U}(1,0,\frac{\lambda_1}{\lambda_2})\quad\text{and}\quad 
			p_j=\frac{\lambda_1}{\lambda_2}\Gamma(j)\mathbf{U}(j,0,\frac{\lambda_1}{\lambda_2}),\quad j\in \mathbb{N}\setminus \{0\}.
		\end{equation}
	\end{enumerate}
\end{teo}

To obtain Theorem \ref{teo1} (i)-(iv) from Theorem \ref{teo2} we can use known asymptotics for $\mathbf{B}$, ${}_2\mathbf{F}_1$ and $\mathbf{U}$. For the convenience of the reader comments are provided in Appendix \ref{S:Uasymp}.

\begin{appendix}
\section{Some known key facts}\label{S:key}

The first fact was observed in \cite[Theorem 1 and its proof]{Neuts}: Given $0\leq s_1\leq ...\leq s_n\leq u$, we have 
\begin{multline}\label{E:Neuts}
	\hat{\mathbb{P}}(\hat{\tau}_1^\ast\leq s_1,...,\hat{\tau}_n^\ast\leq s_n  \mid \ell(\hat Z_u)=n)
	=\hat{\mathbb{P}}(\hat{\tau}_1^\ast\leq s_1,...,\hat{\tau}_n^\ast\leq s_n\mid \hat{\tau}^\ast_n<u\leq \hat{\tau}^\ast_{n+1})\\
	=\frac{\lambda_1^n e^{-\lambda_1 nu}n!}{(1-e^{-\lambda_1u})^kn}\int_0^{s_1}\int_{v_1}^{s_2}\cdots \int_{v_{n-1}}^{s_n}e^{\lambda_1\sum_{i=1}^nv_i}dv_n\cdots dv_2dv_1.
\end{multline}

The next fact was used in \cite{politosacerdotelanski} and the sources quoted there: Suppose that $T_1,...,T_n$ are independent $[0,u]$-valued random variables over some probability space with probability measure $\mathbb{Q}$ and all have the same density
$f_u(s)=\lambda_1e^{\lambda_1 s}\big(e^{\lambda_1 u}-1\big)^{-1},\quad 0\leq s\leq u$.

Now let $T_{(1)},...,T_{(n)}$ denote their order statistics. Then $\mathbb{Q}(T_{(1)}\leq s_1,...,T_{(n)}\leq s_n)$ equals (\ref{E:Neuts}). Consequently, given $0\leq s\leq u$, we have 
\[\frac{1}{n}\sum_{j=1}^n \hat{\mathbb{P}}(\hat{\tau}^\ast_j\leq s\mid \ell(\hat Z_u)=n)=\frac{1}{n}\sum_{j=1}^n \mathbb{Q}(T_{(j)}\leq s)=\mathbb{Q}(T_1\leq s).\]
That is, $f_u$ is the density of the formation time of a household chosen with uniform probability among all the households existing at time $u$. Consequently, using the Markov property of the birth-death process $Z=(Z_u)_{u\geq 0}$ starting at one over $(\Omega',\mathcal{F}',\mathbb{P}')$, 
\begin{multline}\label{LimZ}
	p_j= \lim_{u\rightarrow\infty}\hat{\mathbb{P}}(s(\hat{Z}_u)=j)
	= \lim_{u\rightarrow\infty} \int_{0}^{u}\mathbb{P}'(Z_{u-s}=j \mid Z_s=1) f_u(s)ds\\
	= \lim_{u\rightarrow\infty}  \frac{\lambda_1}{1-e^{-\lambda_1 u}}  \int_{0}^{u} e^{-\lambda_1 v} \mathbb{P}'(Z_{v}=j)dv
	= \lambda_1  \int_{0}^{\infty} e^{-\lambda_1 v} \mathbb{P}'(Z_{v}=j)dv,\quad j\in \mathbb{N}.
\end{multline}	

For	the pure birth case $\mu_2=0$ one can now use the classical fact that 
\[\mathbb{P}'(Z_{v}=0)=0\quad \text{and}  \quad\mathbb{P}'(Z_{v}=j)=e^{-\lambda_2 v}(1-e^{-\lambda_2 v})^{j-1},\ j\in \mathbb{N}\setminus \{0\},\] see for instance \cite[Chapter 8]{Bailey64}. This gives the well-known expression (\ref{nuZero}) in terms of $\mathbf{B}$, cf. \cite[Section 2]{politosacerdotelanski} or \cite{Yule1925}. Now suppose that $\mu_2>0$. In the super- and subcritical cases one can use the known fact that 
\[\mathbb{P}'(Z_{v}=0)=\frac{\mu_2(1-e^{(\lambda_2-\mu_2)u})}{\lambda_2-\mu_2e^{-(\lambda_2-\mu_2)u}}\]
and
\[\mathbb{P}'(Z_{v}=j)=(\lambda_2-\mu_2)^2e^{-(\lambda_2-\mu_2)u}\frac{\lambda_2^{j-1}(1-e^{-(\lambda_2-\mu_2)u})^{j-1}}{(\lambda_2-\mu_2e^{-(\lambda_2-\mu_2)u})^{j+1}},\ j\in \mathbb{N}\setminus \{0\},\]
see \cite{Bailey64}[Chapter 8]. Plugging these expressions into (\ref{LimZ}), one can rewrite the result in terms of ${}_2\mathbf{F}_1$ to obtain (\ref{supercriticalp0}), (\ref{supercriticalpj}), (\ref{subcriticalp0}) and (\ref{subcriticalpj}). This was observed in \cite[Section 3.2]{politosacerdotelanski}. In the critical case one can use 
\[\mathbb{P}'(Z_{v}=0)=\frac{\lambda_2 u}{1+\lambda_2 u}\quad \text{and} \quad \mathbb{P}'(Z_{v}=j)=\frac{(\lambda_2 u)^{j-1}}{(1+\lambda_2 u)^{j+1}},\ j\in \mathbb{N}\setminus \{0\},\]
\cite[Chapter 8]{Bailey64}. Plugging into (\ref{LimZ}), one obtains (\ref{criticalp0}), as found in \cite[Section 3.1]{politosacerdotelanski}. 

To see (\ref{Expectation}), let $\mathbb{E}'$ denote the expectation with respect to $\mathbb{P}'$ and note that Fubini and (\ref{LimZ}) give 
\[\mathbb{E}[d(V_\infty)]=\sum_{j=0}^\infty j p_j=\lambda_1\int_0^\infty e^{-\lambda_1 v}\sum_{j=0}^\infty j\mathbb{P}'(Z_v=j)dv=\lambda_1\int_0^\infty e^{-\lambda_1 v} \mathbb{E}'[Z_v]dv.\]
Now it suffices to take into account the known fact that 
\[\bbE'(Z_v)=
\begin{cases}
	e^{(\lambda_2-\mu_2)v}, & \lambda_2\neq \mu_2\\
	1, & \lambda_2= \mu_2,
\end{cases}\]
see \cite[Chapter 8]{Bailey64}.

\begin{obs} In a similar manner one can calculate the variance $\mathbb{V}[d(V_\infty)]$ of a random variable $d(V_\infty)$ with distribution $(p_j)_{j\geq 0}$, this gives
	\[
	\mathbb{V}[d(V_\infty)]=
	\begin{cases}
		\frac{\frac{2\lambda_1\lambda_2}{\lambda_1-2(\lambda_2-\mu_2)}-\frac{\lambda_1(\lambda_2+\mu_2)}{\lambda_1-(\lambda_2-\mu_2)}}{\lambda_2-\mu_2} - \Big(\frac{\lambda_1}{\lambda_1-(\lambda_2-\mu_2)}\Big)^2, & \text{ if } \mu_2>\lambda_2 \text{ or } \lambda_2-2\lambda_1<\mu_2 <\lambda_2\\<
		2\Big(\frac{\lambda_2}{\lambda_1}\Big), & \text{ if } \lambda_2=\mu_2\\
		\infty, &\text{ otherwise}.
	\end{cases}
	\]
\end{obs}

\section{Large parameter asymptotics}\label{S:Uasymp}

For fixed $b>0$ we have $\mathbf{B}(j,b)\sim \Gamma(b)j^{-b}$ as $j\to \infty$. Combined with Theorem \ref{teo2} (i) this gives Theorem \ref{teo1} (i).

To obtain the asymptotics in Theorem \ref{teo1} (ii) and (iii) from (\ref{supercriticalpj}) and (\ref{subcriticalpj}), respectively, note that for any $b>0$, $\gamma>0$ and $0<z<1$, we have 
\[{}_2\mathbf{F}_1(j+1,b,j+1+\gamma,z)={}_2\mathbf{F}_1(b, j+1,j+1+\gamma,z)\sim (1-z)^{-b},\]
see \cite[Section 2]{Temme03}.

To see Theorem \ref{teo1} (iv) we use the fact that for any $z>0$ there is a function $R(j)$ of $j$ and a constant $C>0$ such that $\limsup_{j\to\infty} R(j)<C$ and
\begin{equation}\label{E:largepar}
	\Gamma(j)\mathbf{U}(j,0,z)=2e^{\frac{z}{2}}\left(\sqrt{2\beta\tanh\Big(\frac{w}{2}\Big)}K_1(2j\beta)+\frac{\frac{1}{j}+\beta}{j+j\beta}e^{-2j\beta}R(j)\right),
\end{equation}
here $w=w(j)$ and $\beta=\beta(j)$ are given by 
$w(j)=\cosh^{-1}\Big(1+\frac{z}{2j}\Big)$ and $ \beta(j)=\frac12(w(j)+\sinh(w(j))),$
and $K_1$ denotes the modified Bessel function of the second kind of order one, see \cite[13.8 (iii)]{DLMF}. Since $w(j)\sim \frac{z}{2j}+\sqrt{\Big(1+\frac{z}{2j}\Big)^2-1}$ and $\beta(j)\sim w(j)$, we find that 
\begin{equation}\label{E:mix}
	2j\beta(j)\sim z+2\sqrt{jz+\frac{z^2}{4}}\sim 2\sqrt{jz}\quad \text{and}\quad \frac{\beta(j)}{j}\sim \frac{\sqrt{z}}{j^{3/2}}.  
\end{equation}
The known asymptotics for $K_1$ gives $K_1(2j\beta(j))\sim \sqrt{\frac{\pi}{4 j\beta(j)}}e^{-2j\beta(j)}$, see \cite[10.25 (ii)]{DLMF}. Consequently the first summand in brackets on the right-hand side of (\ref{E:largepar}) behaves like 
\[\sim \beta(j)K_1(2j\beta(j))\sim \frac{\sqrt{\pi}}{2}\sqrt{\frac{\beta(j)}{j}}e^{-2j\beta(j)}.\]
Given $m\in \mathbb{N}$ we have $j^me^{-2j\beta(j)}<j^m(2j\beta(j))^{-2m}\leq z^{-m}$ for all large enough $j$ by the first part in (\ref{E:mix}). Using also the second part we see that $j^m \beta(j)K_1(2j\beta(j))\to 0$.
For any $\varepsilon>0$ we eventually have $e^{\varepsilon j}e^{-2j\beta(j)}<e^{\frac{\varepsilon}{2}j}$ and therefore $e^{\varepsilon j}\beta(j)K_1(2j\beta(j))\to +\infty$.
The second summand in brackets on the right-hand side of (\ref{E:largepar}) behaves similarly.
\end{appendix}

\end{document}